\documentclass[english,11pt]{article}
\usepackage[latin1]{inputenc} 
\usepackage{amsmath,amsthm,amssymb,babel}

\newtheorem{teo}{Theorem}

\newtheorem{prop}{Proposition}

\newtheorem{lemma}{Lemma}

\def\proof{{\it Proof.}\ }

\def\eq#1{(\ref{#1})}

\def\neweq#1{\begin{equation}\label{#1}}
\def\endeq{\end{equation}}

\def\phi{\varphi}
\def\RR{{\mathbb R} }

\def\di{\displaystyle}
\def\ri{\rightarrow}

\date{}

\title{\sc A multiplicity result for a nonlinear degenerate problem 
arising in the theory of electrorheological fluids\thanks{
Correspondence address: Vicen\c{t}iu R\u{a}dulescu, Department of Mathematics, University of 
Craiova,  200585 Craiova, Romania. E-mail: {\tt radulescu@inf.ucv.ro}}}
\author{\sc Mihai Mih\u ailescu and Vicen\c{t}iu R\u{a}dulescu\\
\small
Department of Mathematics, University of Craiova,  200585 Craiova, Romania\\
\small
E-mail addresses: {\tt mmihailes@yahoo.com}\qquad {\tt radulescu@inf.ucv.ro}}

\begin{document}
\baselineskip16pt 
\maketitle
\noindent{\small{\sc Abstract}.
We study the boundary value problem $-{\rm div}(a(x,\nabla u))=\lambda(u^{\gamma-1}-
u^{\beta-1})$ in $\Omega$,
$u=0$ on $\partial\Omega$, where $\Omega$ is a smooth bounded domain in $\RR^N$ and
${\rm div}(a(x,\nabla u))$ is a $p(x)$-Laplace type operator, with
$1<\beta<\gamma<\inf_{x\in\Omega}p(x)$. We prove that if $\lambda$ is large enough
then there exist at least two nonnegative weak solutions. Our approach relies on the variable 
exponent theory of generalized Lebesgue-Sobolev spaces, combined with adequate variational methods 
and a variant of Mountain Pass Lemma.  \\
\small{\bf 2000 Mathematics 
Subject Classification:}  35D05, 35J60, 35J70, 58E05, 68T40, 76A02. \\ 
\small{\bf Key words:}  $p(x)$-Laplace operator, generalized Lebesgue-Sobolev
space, critical point, weak solution, electrorheological fluids.}

\section{Introduction and preliminary results}
Most materials can be modeled with sufficient accuracy using classical Lebesgue and Sobolev 
spaces, $L^p$ and $W^{1,p}$, where $p$ is a fixed constant. For some materials with 
inhomogeneities, for instance electrorheological fluids (sometimes referred to as ``smart 
fluids"), this is not adequate, but rather the exponent $p$ should be able to vary. This leads us 
to the study of variable exponent Lebesgue and Sobolev spaces, $L^{p(x)}$ and $W^{1,p(x)}$, where 
$p$ is a real--valued function.

This paper is motivated by phenomena which are described by nonlinear boundary value problems of 
the type
\begin{equation}\label{1}
\left\{\begin{array}{lll}
-{\rm div}(a(x,\nabla u))=f(x,u), &\mbox{for}& 
x\in\Omega\\
u=0, &\mbox{for}& x\in\partial\Omega
\end{array}\right.
\end{equation}
where $\Omega\subset\RR^N$ ($N\geq 3$) is a bounded domain with 
smooth boundary, $1<p(x)$ and $p(x)\in C(\overline\Omega)$. 
The interest in studying such problems consists in the presence 
of the $p(x)$--Laplace type operator ${\rm div}(a(x,\nabla u))$.
We remember that the $p(x)$-Laplace operator is defined by
$\Delta_{p(x)}u={\rm div}(|\nabla u|^{p(x)-2}\nabla u)$. 
The study of differential equations
and variational problems involving $p(x)$-growth conditions is
a consequence of their applications. 
Materials requiring such more advanced theory have been studied experimentally since the middle of 
last century. The first major discovery on electrorheological fluids is due to Willis Winslow in  
1949. These fluids have the
interesting property that their viscosity depends on the electric field in the fluid. He noticed 
that such fluids (for instance lithium polymetachrylate) viscosity in an electrical field is 
inversely proportional to the strength of the field. The field induces string-like formations in 
the fluid, which are parallel to the field. They can raise the viscosity by as much as five orders 
of magnitude. This phenomenon is known as the Winslow effect. For a
general account of the underlying physics confer \cite{hal} and for some technical applications
\cite{pfe}.
Electrorheological fluids have been used in robotics and space technology. The experimental 
research has been done mainly in the USA, for instance in NASA laboratories. For more information 
on properties, modelling and the application of variable exponent spaces to these fluids we refer 
to \cite{AM,CF,D,FZZ,hal,R}.

Variable exponent Lebesgue spaces appeared in the literature for the first time already in a 1931 
article by W.~Orlicz \cite{orl}, who proved various results (including H\"older's inequality) in a 
discrete framework.
Orlicz also considered the variable exponent function space $L^{p(x)}$ on the real line, and 
proved the H\"older inequality in this setting, too. Next, Orlicz abandoned the study of variable 
exponent spaces, to concentrate on the theory of the function spaces that now bear his name.
The first systematic study of spaces with variable exponent (called {\it modular spaces}) is due 
to Nakano \cite{nak}. In the appendix of this book, Nakano mentions explicitly variable exponent 
Lebesgue spaces as an example of the more general spaces he considers \cite[p.~284]{nak}.
Despite their broad interest, these spaces have not reached the same main-stream position as 
Orlicz spaces. Somewhat later, a more explicit version of such spaces, namely modular function 
spaces, were investigated by Polish mathematicians. We refer to the book of Musielak \cite{M} for 
a nice presentation of the modular function spaces. This book, although not dealing specifically 
with the spaces that interest us, is still specific enough to contain several interesting results 
regarding variable exponent spaces. 
Variable exponent Lebesgue spaces on the real line have been independently developed by Russian 
researchers, notably Sharapudinov. These investigations originated in a paper by Tsenov   
\cite{tse}. The question raised by Tsenov and solved by Sharapudinov \cite{sha} is the 
minimization of
$\int_a^b \vert u(x)-v(x)\vert^{p(x)}dx, $
where $u$ is a fixed function and $v$ varies over a finite dimensional subspace of 
$L^{p(x)}([a,b])$.  Sharapudinov also introduces the Luxemburg norm for the Lebesgue space and 
shows that this space is reflexive if the exponent satisfies $1 < p^- \leq p^+ < \infty$. In the 
80's Zhikov started a new line of investigation, that was to become intimately related to the 
study of variable exponent spaces, namely he considered variational integrals with non-standard 
growth conditions. 

We recall in what follows some definitions and basic properties
of the generalized Lebesgue--Sobolev spaces $L^{p(x)}(\Omega)$
and $W_0^{1,p(x)}(\Omega)$, where $\Omega$ is
a bounded domain in $\RR^N$. 

Throughout this paper we assume that $p(x)>1$, $p(x)\in C^{0,\alpha}
(\overline\Omega)$ with $\alpha\in(0,1)$.

Set
$$C_+(\overline\Omega)=\{h;\;h\in C(\overline\Omega),\;h(x)>1\;{\rm for}\;
{\rm all}\;x\in\overline\Omega\}.$$
For any $h\in C_+(\overline\Omega)$ we define
$$h^+=\sup_{x\in\Omega}h(x)\qquad\mbox{and}\qquad h^-=
\inf_{x\in\Omega}h(x).$$
For any $p(x)\in C_+(\overline\Omega)$, we define the variable exponent Lebesgue space
$$L^{p(x)}(\Omega)=\{u;\ u\ \mbox{is a
 measurable real-valued function such that }
\int_\Omega|u(x)|^{p(x)}\;dx<\infty\}.$$
We define a norm, the so-called {\it Luxemburg norm}, on this space by the formula
$$|u|_{p(x)}=\inf\left\{\mu>0;\;\int_\Omega\left|
\frac{u(x)}{\mu}\right|^{p(x)}\;dx\leq 1\right\}.$$    
Variable exponent Lebesgue spaces resemble classical Lebesgue spaces
in many respects: they are Banach spaces \cite[Theorem 2.5]{KR}, the H\"older inequality holds
\cite[Theorem 2.1]{KR}, they are reflexive if and only if $1 < p^-\leq p^+<\infty$  
\cite[Corollary 2.7]{KR} and continuous functions are dense if $p^+ <\infty$ \cite[Theorem 
2.11]{KR}. The inclusion between
Lebesgue spaces also generalizes naturally \cite[Theorem 2.8]{KR}: if $0 < |\Omega|<\infty$
 and $p$, $q$
are variable exponent so that $p_1(x) \leq p_2(x)$ almost everywhere in $\Omega$ then there exists 
the continuous embedding
$L^{p_2(x)}(\Omega)\hookrightarrow L^{p_1(x)}(\Omega)$, whose norm does not exceed $|\Omega|+1$.

We denote by $L^{q(x)}(\Omega)$ the conjugate space 
of $L^{p(x)}(\Omega)$, where $1/p(x)+1/q(x)=1$. For any 
$u\in L^{p(x)}(\Omega)$ and $v\in L^{q(x)}(\Omega)$ the H\"older 
type inequality
\begin{equation}\label{Hol}
\left|\int_\Omega uv\;dx\right|\leq\left(\frac{1}{p^-}+
\frac{1}{q^-}\right)|u|_{p(x)}|v|_{q(x)}
\end{equation}
holds true.  

An important role in manipulating the generalized Lebesgue-Sobolev spaces is played by the {\it 
modular} of the $L^{p(x)}(\Omega)$ space, which is the mapping 
 $\rho_{p(x)}:L^{p(x)}(\Omega)\rightarrow\RR$ defined by
$$\rho_{p(x)}(u)=\int_\Omega|u|^{p(x)}\;dx.$$
If $(u_n)$, $u\in L^{p(x)}(\Omega)$ and $p^+<\infty$ then the following relations 
holds true
\begin{equation}\label{L4}
|u|_{p(x)}>1\;\;\;\Rightarrow\;\;\;|u|_{p(x)}^{p^-}\leq\rho_{p(x)}(u)
\leq|u|_{p(x)}^{p^+}
\end{equation}  
\begin{equation}\label{L5}
|u|_{p(x)}<1\;\;\;\Rightarrow\;\;\;|u|_{p(x)}^{p^+}\leq
\rho_{p(x)}(u)\leq|u|_{p(x)}^{p^-}
\end{equation}
\begin{equation}\label{L6}
|u_n-u|_{p(x)}\rightarrow 0\;\;\;\Leftrightarrow\;\;\;\rho_{p(x)}
(u_n-u)\rightarrow 0.
\end{equation}
Spaces with $p^+ =\infty$ have been studied by Edmunds, Lang and Nekvinda \cite{edm}.

Next, we define $W_0^{1,p(x)}(\Omega)$ as the closure of 
$C_0^\infty(\Omega)$ under the norm 
$$\|u\|=|\nabla u|_{p(x)}.$$
The space $(W_0^{1,p(x)}(\Omega),\|\cdot\|)$ is a separable and 
reflexive Banach space. We note that if $q\in C_+(\overline\Omega)$ 
and $q(x)<p^\star(x)$ for all $x\in\overline\Omega$ then the 
embedding
$W_0^{1,p(x)}(\Omega)\hookrightarrow L^{q(x)}(\Omega)$ 
is compact and continuous, where $p^\star(x)=\frac{Np(x)}{N-p(x)}$ 
if $p(x)<N$ or $p^\star(x)=+\infty$ if $p(x)\geq N$. We refer to \cite{edm2,edm3,FZ1,KR} for 
further properties of variable exponent Lebesgue-Sobolev spaces.
 
 \section{The main result}
Assume that $a(x,\xi):\overline\Omega\times\RR^N\rightarrow\RR^N$ 
is the continuous derivative with respect to $\xi$ of 
the mapping $A:\overline\Omega\times\RR^N\rightarrow\RR$, 
$A=A(x,\xi)$, that is, $a(x,\xi)=\frac{d}{d\xi}A(x,\xi)$. Suppose that 
$a$ and $A$ satisfy the following hypotheses:
\smallskip

\noindent (A1) The following equality holds 
$$A(x,0)=0,$$ 
for all $x\in\overline\Omega$.
\smallskip

\noindent (A2) There exists a positive constant $c_1$ such that 
$$|a(x,\xi)|\leq c_1(1+|\xi|^{p(x)-1}),$$ 
for all $x\in\overline\Omega$ and $\xi\in\RR^N$.
\smallskip

\noindent (A3) The following inequality holds
$$0\leq(a(x,\xi)-a(x,\psi))\cdot(\xi-\psi),$$
for all $x\in\overline\Omega$ and $\xi,\psi\in\RR^N$, with equality
if and only if $\xi=\psi$.

\noindent (A4) There exists $k>0$ such that
$$A\left(x,\frac{\xi+\psi}{2}\right)\leq\frac{1}{2}A(x,\xi)+
\frac{1}{2}A(x,\psi)-k|\xi-\psi|^{p(x)}$$
for all $x\in\overline\Omega$ and $\xi,\psi\in\RR^N$.
\smallskip

\noindent (A5) The following inequalities hold true 
$$|\xi|^{p(x)}\leq a(x,\xi)\cdot\xi\leq p(x)\;A(x,\xi),$$ 
for all $x\in\overline\Omega$ and $\xi\in\RR^N$.
\bigskip

\noindent{\bf Examples.}\\
1. Set $A(x,\xi)=\frac{1}{p(x)}|\xi|^{p(x)}$, $a(x,\xi)=|\xi|^{p(x)-2}
\xi$, where $p(x)\geq 2$. Then we get the $p(x)$-Laplace operator
$${\rm div}(|\nabla u|^{p(x)-2}\nabla u).$$
\smallskip

\noindent 2. Set $A(x,\xi)=\frac{1}{p(x)}[(1+|\xi|^2)^{p(x)/2}-1]$,
$a(x,\xi)=(1+|\xi|^2)^{(p(x)-2)/2}\xi$, where
$p(x)\geq 2$. Then we obtain the generalized mean curvature operator
$${\rm div}((1+|\nabla u|^2)^{(p(x)-2)/2}\nabla u).$$
\bigskip

In this paper we study problem \eq{1} in the particular case 
$$f(x,t)=\lambda(t^{\gamma-1}-t^{\beta-1})$$ 
with $1<\beta<\gamma<\inf_{x\in\overline\Omega}p(x)$ and $t\geq 0$. 
More precisely, we consider the degenerate boundary value problem 
\begin{equation}\label{2}
\left\{\begin{array}{lll}
-{\rm div}(a(x,\nabla u))=\lambda(u^{\gamma-1}-
u^{\beta-1}), &\mbox{for}& x\in\Omega\\
u=0, &\mbox{for}& x\in\partial\Omega\\
u\geq 0,&\mbox{for}& x\in\Omega.
\end{array}\right.
\end{equation}

We say that $u\in W_0^{1,p(x)}(\Omega)$ is a {\it weak 
solution} of problem \eq{2} if
$u\geq 0$ a. e. in $\Omega$
and
$$\int_\Omega a(x,\nabla u)\cdot\nabla\phi\;dx-\lambda
\int_\Omega u^{\gamma -1}\phi\;dx+\lambda\int_\Omega 
u^{\beta-1}\phi\;dx=0$$
for all $\phi\in W_0^{1,p(x)}(\Omega)$.
\smallskip

Our main result asserts that problem \eq{2} has at least two nontrivial 
weak solutions provided that $\lambda>0$ is large enough and 
operators $A$ and $a$ satisfy conditions (A1)-(A5). More precisely, 
we prove

\begin{teo}\label{t1}
Assume hypotheses (A1)-(A5) are fulfilled. Then there exists 
$\lambda^\star>0$ such that for all $\lambda>\lambda^\star$ 
problem \eq{2} has at least two distinct non-negative,
nontrivial weak solutions, provided that $p^+<\min\{N, 
Np^-/(N-p^-)\}$. 
\end{teo}

{\bf Remark.} By Theorem 4.3 in \cite{FZh} problem \eq{2} has 
at least a weak solution in the particular case 
$a(x,\xi)=|\xi|^{p(x)-1}\xi$. However, the proof in \cite{FZh}
does not state the fact that the solution is non-negative and not even
nontrivial in the case when $f(x,0)=0$.
\smallskip

We point out that our result is inspired 
by  \cite[Theorem 1.2]{P}, where a related property 
 is proved in the case of the $p$-Laplace operators. We
point out that the extension from $p$-Laplace operator to 
$p(x)$-Laplace operator is not trivial, since the $p(x)$-Laplacian
has a more complicated structure than the $p$-Laplace operator,
for example it is inhomogeneous.

\section{Proof of Theorem \ref{t1}}
Let $E$ denote the generalized Sobolev 
space $W_0^{1,p(x)}(\Omega)$.

Define the energy functional
$I:E\rightarrow\RR$ by
$$I(u)=\int_\Omega A(x,\nabla u)\;dx-\frac{\lambda}
{\gamma}\int_\Omega u_+^\gamma\;dx+\frac{\lambda}{\beta}\int_\Omega
u_+^\beta\;dx\,,$$
where $u_+(x)=\max\{u(x),0\}$.

We first establish some basic properties of $I$.

\begin{prop}\label{p1}
The functional $I$ is well-defined on $E$ and $I\in C^1(E,\RR)$ with
the derivative given by
$$\langle I^{'}(u),\phi\rangle=\int_\Omega a(x,\nabla u)
\cdot\nabla\phi\;dx-\lambda\int_\Omega u_+^{\gamma-1}\phi\;dx
+\lambda\int_\Omega u_+^{\beta-1}\phi\;dx,$$
for all $u$, $\phi\in E$.
\end{prop} 
With that end in view we define the functional $\Lambda:E\rightarrow
\RR$ by
$$\Lambda(u)=\int_{\Omega}A(x,\nabla u)\;dx,\;\;\;\forall u\in E.$$
\begin{lemma}\label{l1}
(i) The functional $\Lambda$ is well-defined on $E$.\\
(ii) The functional $\Lambda$ is of class $C^1(E,\RR)$ and
$$\langle\Lambda^{'}(u),\phi\rangle=\int_\Omega a(x,\nabla u)\cdot
\nabla\phi\;dx,$$
for all $u,\phi\in E$.
\end{lemma}
\proof
(i) For any $x\in\Omega$ and $\xi\in\RR^N$ we have
$$A(x,\xi)=\int_0^1\frac{d}{dt}A(x,t\xi)\;dt=\int_0^1a(x,t\xi)\cdot
\xi\;dt.$$
Using hypotheses (A2) we get
\begin{equation}\label{ecc4}
\begin{array}{lll}
A(x,\xi)&\leq&c_1\di\int_0^1(1+|\xi|^{p(x)-1}t^{p(x)-1})|\xi|\;dt\\
&\leq&c_1|\xi|+\di\frac{c_1}{p(x)}|\xi|^{p(x)}\\
&\leq&c_1|\xi|+\di\frac{c_1}{p^-}|\xi|^{p(x)},\;\;\;\forall x\in
\overline\Omega,\;\xi\in\RR^N.
\end{array}
\end{equation}
The above inequality and (A5) imply
$$0\leq\int_\Omega A(x,\nabla u)\;dx\leq c_1\int_\Omega|\nabla u|
\;dx+\frac{c_1}{p^-}\int_\Omega|\nabla u|^{p(x)}\;dx,\;\;\;
\forall u\in E.$$
Using inequality \eq{Hol} and relations \eq{L4} and \eq{L5} we 
deduce that $\Lambda$ is well defined on $E$.
\smallskip

\noindent (ii) {\bf Existence of the G\^ateaux derivative.} Let
$u$, $\phi\in E$. Fix $x\in\Omega$ and $0<|r|<1$. Then, by the mean value
theorem, there exists $\nu\in[0,1]$ such that
$$|A(x,\nabla u(x)+r\nabla\phi(x))-A(x,\nabla u)|/|r|=
|a(x,\nabla u(x)+\nu\;r\;\nabla\phi(x))||\nabla\phi(x)|.$$
Using condition (A2) we obtain
\begin{eqnarray*}
|A(x,\nabla u(x)+r\nabla\phi(x))-A(x,\nabla u)|/|r|&\leq&
[c_1+c_1(|\nabla u(u)|+|\nabla\phi(x)|)^{p(x)-1}]|\nabla\phi(x)|\\
&\leq&[c_1+c_12^{p^+}(|\nabla u(x)|^{p(x)-1}+|\nabla\phi(x)|^{p(x)
-1})]|\nabla\phi(x)|.
\end{eqnarray*}
Next, by inequality \eq{Hol}, we have
$$\int_\Omega c_1|\nabla\phi|\;dx\leq|c_1|_{\frac{p(x)}{p(x)-1}}
\cdot|\nabla\phi|_{p(x)}$$
and
$$\int_\Omega|\nabla u|^{p(x)-1}|\nabla\phi|\;dx\leq
||\nabla u|^{p(x)-1}|_{\frac{p(x)}{p(x)-1}}\cdot|\nabla\phi|_{p(x)}.$$
The above inequalities imply
$$c_1[1+2^{p^+}(|\nabla u(x)|^{p(x)-1}+|\nabla\phi(x)|^{p(x)-1})]
|\nabla\phi(x)|\in L^1(\Omega).$$
It follows from the Lebesgue theorem that
$$\langle\Lambda^{'}(u),\phi\rangle=\int_\Omega a(x,\nabla u)\cdot
\nabla\phi\;dx.$$
{\bf Continuity of the G\^ateaux derivative.} Assume $u_n\rightarrow u$
in $E$. Let us define $\theta(x,u)=a(x,\nabla u)$. Using hypotheses
(A2) and Proposition 2.2 in \cite{FZh} we deduce that $\theta(x,u_n)
\rightarrow\theta(x,u)$ in $(L^{q(x)}(\Omega))^N$, where $q(x)=
\frac{p(x)}{p(x)-1}$. By inequality \eq{Hol} we obtain
$$|\langle\Lambda^{'}(u_n)-\Lambda^{'}(u),\phi\rangle|\leq
|\theta(x,u_n)-\theta(x,u)|_{q(x)}|\nabla\phi|_{p(x)}$$
and so
$$\|\Lambda^{'}(u_n)-\Lambda^{'}(u)\|\leq|\theta(x,u_n)-
\theta(x,u)|_{q(x)}\rightarrow 0,\;\;\;{\rm as}\;n\rightarrow
\infty.$$
The proof of Lemma \ref{l1} is complete.  \qed

\begin{lemma}\label{l2}
If $u\in E$ then $u_+$, $u_-\in E$ and
$$\nabla u_+=\left\{\begin{array}{lll}
0, &\mbox{if}& [u\leq 0]\\
\nabla u, &\mbox{if}& [u>0],
\end{array}\right. \qquad
\nabla u_-=\left\{\begin{array}{lll}
0, &\mbox{if}& [u\geq 0]\\
\nabla u, &\mbox{if}& [u<0]
\end{array}\right.$$
where $u_\pm=\max\{\pm u(x),0\}$ for all $x\in\Omega$.
\end{lemma}

\proof
Let $u\in E$ be fixed. Then there exists a sequence $(\phi_n)\in 
C_0^\infty(\Omega)$ such that
$$|\nabla(\phi_n-u)|_{p(x)}\rightarrow 0.$$
Since $1<p^-\leq p(x)$ for all $x\in\Omega$, it follows  that 
$L^{p(x)}$ is continuously embedded in $L^{p^-}(\Omega)$ and thus
$$|\nabla(\phi_n-u)|_{p^-}\rightarrow 0.$$
Hence $u\in W^{1,p^-}_0(\Omega)$. We obtain
\begin{equation}\label{ecca}
u_+,u_-\in W_0^{1,p^-}(\Omega)\subset W_0^{1,1}(\Omega).
\end{equation}
On the other hand, Theorem 7.6 in \cite{GT} implies
$$\nabla u_+=\left\{\begin{array}{lll}
0, &\mbox{if}& [u\leq 0]\\
\nabla u, &\mbox{if}& [u>0],
\end{array}\right.
\nabla u_-=\left\{\begin{array}{lll}
0, &\mbox{if}& [u\geq 0]\\
\nabla u, &\mbox{if}& [u<0].
\end{array}\right.$$
By the above equalities we deduce that
\begin{equation}\label{ecc9}
|u_+(x)|^{p(x)}\leq|u(x)|^{p(x)},\;\;\;
|\nabla u_+(x)|^{p(x)}\leq|\nabla u|^{p(x)},\;\;\;{\rm a.e.}\;
x\in\Omega
\end{equation}
and
\begin{equation}\label{ecc10}
|u_-(x)|^{p(x)}\leq|u(x)|^{p(x)},\;\;\;
|\nabla u_-(x)|^{p(x)}\leq|\nabla u|^{p(x)},\;\;\;{\rm a.e.}\;
x\in\Omega.
\end{equation}
Since $u\in E$ we have
\begin{equation}\label{ecc11}
|u(x)|^{p(x)},\;|\nabla u(x)|^{p(x)}\in L^1(\Omega).
\end{equation}
By \eq{ecc9}, \eq{ecc10} and \eq{ecc11} and Lebesgue theorem we 
obtain that $u_+$, $u_-\in L^{p(x)}(\Omega)$ and $\rho_{p(x)}
(|\nabla u_+|)<\infty$, $\rho_{p(x)}(|\nabla u_-|)<\infty$. 
It follows that
\begin{equation}\label{eccb}
u_+,u_-\in W^{1,p(x)}(\Omega)
\end{equation}  
where $W^{1,p(x)}(\Omega)=\{u\in L^{p(x)}(\Omega);\;\;|\nabla u|
\in L^{p(x)}(\Omega)\}$ (see \cite{FZ1} for more details).

By \eq{ecca} and \eq{eccb} we conclude that
$$u_+,u_-\in W^{1,p(x)}(\Omega)\cap W_0^{1,1}(\Omega).$$
Since $p\in C^{0,\alpha}(\overline\Omega)$,
Theorem 2.6 and Remark 2.9 in \cite{FZ1} show that $E=
W^{1,p(x)}(\Omega)\cap W_0^{1,1}(\Omega)$. Thus 
$u_+$, $u_-\in E$ and
the proof of Lemma \ref{l2} is complete.  \qed

By Lemmas \ref{l1} and \ref{l2} it is clear that Proposition \ref{p1}
holds true. Thus, the weak solutions of \eq{2} are exactly the 
critical points of $I$. The above remark shows that we can prove 
Theorem \ref{t1} using the critical points theory. More exactly, we 
first show that for $\lambda>0$ large enough, the functional $I$ has 
a global minimizer $u_1\geq 0$ such that $I(u_1)<0$. Next,
by means of the Mountain 
Pass Theorem, a second 
critical point $u_2$ with $I(u_2)>0$ is obtained.
\smallskip

\noindent{\bf Remark.} If $u$ is a critical point of $I$  then 
using Lemma \ref{l2} and condition (A5) we have
\begin{eqnarray*}
0&=&\langle I^{'}(u),u_-\rangle=\int_\Omega a(x,\nabla u)
\cdot\nabla u_-\;dx-\lambda\int_\Omega(u_+)^{\gamma-1}u_-\;dx
+\lambda\int_\Omega(u_+)^{\beta-1}u_-\;dx\\
&=&\int_\Omega a(x,\nabla u)\cdot\nabla u_-\;dx
=\int_\Omega a(x,\nabla u_-)\cdot\nabla u_-\;dx
\geq\int_\Omega|\nabla u_-|^{p(x)}\;dx.
\end{eqnarray*}
Thus we deduce that $u\geq 0$. It follows that the nontrivial critical
points of $I$ are non-negative solutions of \eq{2}.

\begin{lemma}\label{l3}
The functional $\Lambda$ is weakly lower semi-continuous.
\end{lemma}

\proof
By Corollary III.8 in \cite{B}, it is enough to show that $\Lambda$ 
is  inferior semi-continuous. For this purpose, we fix
$u\in E$ and $\epsilon>0$. Since $\Lambda$ is convex
(by condition (A4)), we deduce that for any $v\in E$ 
the following inequality holds
$$\int_\Omega A(x,\nabla v)\;dx\geq\int_\Omega A(x,\nabla u)\;dx
+\int_\Omega a(x,\nabla u)\cdot(\nabla v-\nabla u)\;dx.$$
Using condition (A2) and inequality \eq{Hol} we have
\begin{eqnarray*}
\int_\Omega A(x,\nabla v)\;dx&\geq&\int_\Omega A(x,\nabla u)\;dx-
\int_\Omega|a(x,\nabla u)||\nabla v-\nabla u|\;dx\\
&\geq&\int_\Omega A(x,\nabla u)\;dx- c_1\int_\Omega|\nabla(v-u)|\;dx
-c_1\int_\Omega|\nabla u|^{p(x)-1}|\nabla(v-u)|\;dx\\
&\geq&\int_\Omega A(x,\nabla u)\;dx-c_2|1|_{q(x)}|\nabla(v-u)|_{p(x)}
-c_3||\nabla u|^{p(x)-1}|_{q(x)}|\nabla(v-u)|_{p(x)}\\
&\geq&\int_\Omega A(x,\nabla u)\;dx-c_4\|v-u\|\\
&\geq&\int_\Omega A(x,\nabla u)\;dx-\epsilon
\end{eqnarray*}
for all $v\in E$ with $\|v-u\|<\delta=\epsilon/c_4$, where $c_2$, $c_3$, $c_4$ are positive 
constants, and 
$q(x)=\frac{p(x)}{p(x)-1}$.
We conclude that $\Lambda$ is weakly lower semi-continuous. The proof 
of Lemma \ref{l3} is complete.  \qed

\begin{lemma}\label{le4}
There exists $\lambda_1>0$ such that
$$\lambda_1=\inf\limits_{u\in E,\;\|u\|>1}\frac{\di\int_\Omega
\di\frac{1}{p(x)}|\nabla u|^{p(x)}\;dx}{\di\int_\Omega|u|^{p^-}\;dx}
.$$
\end{lemma} 
\proof
We know that $E$ is continuously embedded in $L^{p^-}(\Omega)$. It
follows that there exists $C>0$ such that
$$\|u\|\geq C|u|_{p^-},\;\;\;\forall u\in E.$$
On the other hand, by \eq{L4} we have
$$\int_\Omega|\nabla u|^{p(x)}\;dx\geq\|u\|^{p^-},\;\;\;\forall 
u\in E\;{\rm with}\;\|u\|>1.$$
Combining the above inequalities we obtain
$$\int_\Omega\frac{1}{p(x)}|\nabla u|^{p(x)}\;dx\geq\frac{C^{p^-}}
{p^+}\int_\Omega|u|^{p^-}\;dx,\;\;\;\forall u\in E\;{\rm with}\;
\|u\|>1.$$
The proof of Lemma \ref{le4} is complete.\qed

\begin{prop}\label{p2}
(i) The functional $I$ is bounded from below and coercive.\\
(ii) The functional $I$ is weakly lower semi-continuous.
\end{prop}
\proof
(i) Since $1<\beta<\gamma<p^-$ we have
$$\lim\limits_{t\rightarrow\infty}\frac{\di\frac{1}{\gamma}t^\gamma-
\di\frac{1}{\beta}t^\beta}{t^{p^-}}=0.$$ 
Then for any $\lambda>0$ there exists $C_\lambda>0$ such that
$$\lambda\left(\frac{1}{\gamma}t^\gamma-\frac{1}{\beta}t^\beta
\right)\leq\frac{\lambda_1}{2}t^{p^-}+C_\lambda,\;\;\;\forall 
t\geq 0,$$
where $\lambda_1$ is defined in Lemma \ref{le4}.

Condition (A5) and the above inequality show that for any $u\in E$ 
with $\|u\|>1$ we have
\begin{eqnarray*}
I(u)&\geq&\int_\Omega\frac{1}{p(x)}|\nabla u|^{p(x)}\;dx-
\frac{\lambda_1}{2}\int_\Omega|u|^{p^-}\;dx-C_\lambda\mu(\Omega)\\
&\geq&\frac{1}{2}\int_\Omega\frac{1}{p(x)}|\nabla u|^{p(x)}\;dx-
C_\lambda\mu(\Omega)\\
&\geq&\frac{1}{2p^+}\|u\|^{p^-}-C_\lambda\mu(\Omega).
\end{eqnarray*}
This shows that $I$ is bounded from below and coercive.
\smallskip

\noindent (ii) Using Lemma \ref{l3} we deduce that
$\Lambda$ is weakly lower semi-continuous.
We show that $I$ is weakly lower semi-continuous. Let 
$(u_n)\subset E$ be a sequence which converges weakly to $u$ in $E$.
Since $\Lambda$ is weakly lower semi-continuous we have
\begin{equation}\label{3}
\Lambda(u)\leq\liminf\limits_{n\rightarrow\infty}\Lambda(u_n).
\end{equation}
On the other hand, since $E$ is compactly embedded in 
$L^\gamma(\Omega)$ and $L^\beta(\Omega)$ it follows that $({u_n}_+)$
converges strongly to $u_+$ both in $L^\gamma(\Omega)$ and in 
$L^\beta(\Omega)$. This fact together with relation \eq{3} imply
$$I(u)\leq\liminf\limits_{n\rightarrow\infty}I(u_n).$$
Therefore, $I$ is weakly lower semi-continuous.
The proof of Proposition \ref{p2} is complete.  \qed

By Proposition \ref{p2} and Theorem 1.2 in \cite{S} we deduce that 
there exists $u_1\in E$ a global minimizer of $I$. The following result implies that
$u_1\neq 0$, provided that $\lambda$ is sufficiently large.
 
\begin{prop}\label{p3}
There exists $\lambda^\star>0$ such that $\inf_E I<0$.
\end{prop}

\proof
Let $\Omega_1\subset\Omega$ be a compact subset, large enough and
$u_0\in E$ be such that $u_0(x)=t_0$ in $\Omega_1$ and $0\leq u_0(x)
\leq t_0$ in $\Omega\setminus\Omega_1$, where $t_0>1$ is chosen 
such that
$$\frac{1}{\gamma}t_0^\gamma-\frac{1}{\beta}t_0^\beta>0.$$
We have
\begin{eqnarray*}
\frac{1}{\gamma}\int_\Omega u_0^\gamma\;dx-\frac{1}{\beta}\int_\Omega
u_0^\beta\;dx&\geq&\frac{1}{\gamma}\int_{\Omega_1}u_0^\gamma\;dx-
\frac{1}{\beta}\int_{\Omega_1}u_0^\beta\;dx-\frac{1}{\beta}
\int_{\Omega\setminus\Omega_1}u_0^\beta\;dx\\
&\geq&\frac{1}{\gamma}
\int_{\Omega_1}u_0^\gamma\;dx-\frac{1}{\beta}\int_{\Omega_1}u_0^\beta
\;dx-\frac{1}{\beta}\;t_0^\beta\;\mu(\Omega\setminus\Omega_1)>0
\end{eqnarray*}
and thus $I(u_0)<0$ for $\lambda>0$ large enough.
The proof of Proposition \ref{p3} is complete.  \qed

Since Proposition \ref{p3} holds true it follows that $u_1\in E$
is a nontrivial weak solution of problem \eq{2}.
\bigskip

Fix $\lambda\geq\lambda^\star$. Set
$$
g(x,t)=\left\{\begin{array}{lll}
0, &\mbox{for}& t<0\\
t^{\gamma-1}-t^{\beta-1}, &\mbox{for}& 0\leq t\leq u_1(x)\\
u_1(x)^{\gamma-1}-u_1(x)^{\beta-1}, &\mbox{for}& t>u_1(x)
\end{array}\right.
$$
and
$$G(x,t)=\int_0^tg(x,s)\;ds.$$
Define the functional $J:E\rightarrow\RR$ by
$$J(u)=\int_\Omega A(x,\nabla u)\;dx-\lambda \int_\Omega G(x,u)
\;dx.$$
The same arguments as those used for functional $I$ imply that
$J\in C^1(E,\RR)$ and
$$\langle J^{'}(u),\phi\rangle=\int_\Omega a(x,\nabla u)
\cdot\nabla\phi\;dx-\lambda\int_\Omega g(x,u)\phi\;dx,$$
for all $u$, $\phi\in E$.

On the other hand, we point out that if $u\in E$ is a 
critical point of $J$ then $u\geq 0$. The proof can be carried out
as in the case of functional $I$.

Next, we prove
\begin{lemma}\label{l5}
If $u$ is a critical point of $J$ then $u\leq u_1$.
\end{lemma} 
\proof
We have
\begin{eqnarray*}
0&=&\langle J^{'}(u)-I^{'}(u_1),(u-u_1)_+\rangle\\
&=&\int_\Omega(a(x,\nabla u)-a(x,\nabla u_1))
\cdot\nabla(u-u_1)_+\;dx-\lambda\int_\Omega[g(x,u)-(u_1^
{\gamma-1}-u_1^{\beta-1})](u-u_1)_+\;dx\\
&=&\int_{[u>u_1]}(a(x,\nabla u)-a(x,\nabla u_1))
\cdot\nabla(u-u_1)\;dx.
\end{eqnarray*}
By condition (A3) we deduce that the above equality holds if and 
only if $\nabla u=\nabla u_1$. It follows that $\nabla u(x)=
\nabla u_1(x)$ for all $x\in\omega:=\{y\in\Omega;\;\;u(y)>u_1(y)\}$. 
Hence
$$\int_\omega|\nabla(u-u_1)|^{p(x)}\;dx=0$$
and thus
$$\int_\Omega|\nabla(u-u_1)_+|^{p(x)}\;dx=0.$$
By relation \eq{L5} we obtain
$$\|(u-u_1)_+\|=0.$$
Since $u-u_1\in E$ by Lemma \ref{l2} we have that $(u-u_1)_+\in E$.
Thus we obtain that $(u-u_1)_+=0$ in $\Omega$, that is, $u\leq u_1$ in
$\Omega$.
The proof of Lemma \ref{l5} is complete.  \qed

In the following we determine a critical point $u_2\in E$ of $J$
such that $J(u_2)>0$ via the Mountain Pass Theorem. By the above lemma
we will deduce that $0\leq u_2\leq u_1$ in $\Omega$. Therefore 
$$g(x,u_2)=u_2^{\gamma-1}-u_2^{\beta-1}\;\;\;{\rm and}\;\;\;
G(x,u_2)=\frac{1}{\gamma}u_2^\gamma-\frac{1}{\beta}u_2^\beta$$
and thus
$$J(u_2)=I(u_2)\;\;\;{\rm and}\;\;\;J^{'}(u_2)=I^{'}(u_2).$$
More exactly we find
$$I(u_2)>0=I(0)>I(u_1)\qquad\mbox{and}\qquad I^{'}(u_2)=0\,.$$
This shows that $u_2$ is a weak solution of problem \eq{2} such
that $0\leq u_2\leq u_1$, $u_2\neq 0$ and $u_2\neq u_1$.

In order to find $u_2$ described above we prove
\begin{lemma}\label{l6}
There exists $\rho\in(0,\|u_1\|)$ and $a>0$ such that
$J(u)\geq a$, for all $u\in E$ with $\|u\|=\rho.$
\end{lemma}

\proof
Let $u\in E$ be fixed, such that $\|u\|<1$. It is clear that there 
exists $\delta>1$ such that
$$\frac{1}{\gamma}t^\gamma-\frac{1}{\beta}t^\beta\leq 0,\;\;\;
\forall t\in[0,\delta].$$
For $\delta$ given above we define
$$\Omega_u=\{x\in\Omega;\;u(x)>\delta\}.$$
If $x\in\Omega\setminus\Omega_u$ with $u(x)<u_1(x)$ we have
$$G(x,u)=\frac{1}{\gamma}u_+^\gamma-\frac{1}{\beta}u_+^\beta\leq 0.$$
If $x\in\Omega\setminus\Omega_u$ with $u(x)>u_1(x)$ then
$u_1(x)\leq\delta$ and we have
$$G(x,u)=\frac{1}{\gamma}u_1^\gamma-\frac{1}{\beta}u_1^\beta\leq 0.$$
Thus we deduce that
$$G(x,u)\leq 0,\;\;\;{\rm on}\;\Omega\setminus\Omega_u.$$
Provided that $\|u\|<1$ by condition (A5) and relation \eq{L5} we get
\begin{equation}\label{s1}
\begin{array}{lll}
J(u)&\geq&\di\int_\Omega\di\frac{1}{p(x)}|\nabla u|^{p(x)}\;dx-
\lambda\di\int_{\Omega_u}G(x,u)\;dx\\
&\geq&\di\frac{1}{p^+}\|u\|^{p^+}-\lambda\di\int_{\Omega_u}G(x,u)\;dx
\end{array}
\end{equation}
Since $p^+<\min\{N, \frac{Np^-}{N-p^-}\}$ it follows that 
$p^+<p^\star(x)$ for all $x\in\overline\Omega$. Then there exists
$q\in(p^+,\frac{Np^-}{N-p^-})$ such that $E$ is continuously
embedded in $L^q(\Omega)$. Thus there exists a positive constant
$C>0$ such that
$$|u|_q\leq C\|u\|,\;\;\;\forall u\in E.$$
Using the definition of $G$, H\"older's inequality and the 
above estimate, we obtain
\begin{equation}\label{s2}
\begin{array}{lll}
\lambda\di\int_{\Omega_u} G(x,u)\;dx&=&\lambda\di\int_{\Omega_u\cap
[u<u_1]}\left(\di\frac{1}{\gamma}u_+^\gamma-\di\frac{1}{\beta}
u_+^\beta\right)\;dx+\lambda\di\int_{\Omega_u\cap[u>u_1]}\left(
\di\frac{1}{\gamma}u_1^\gamma-\di\frac{1}{\beta}u_1^\beta\right)
\;dx\\
&\leq&\di\frac{2\lambda}{\gamma}\di\int_{\Omega_u}u_+^\gamma\;dx\\
&\leq&\di\frac{2\lambda}{\gamma}\di\int_{\Omega_u}u_+^{p^+}\;dx\\
&\leq&\di\frac{2\lambda}{\gamma}\left(\di\int_{\Omega_u}u_+^{q}
\;dx\right)^{p^+/q}[\mu(\Omega_u)]^{1-p^+/q}\\
&\leq&C\di\frac{2\lambda}{\gamma}[\mu(\Omega_u)]^{1-p^+/q}\|u\|^{p^+}.
\end{array}
\end{equation}
By \eq{s1} and \eq{s2} we infer that it is enough to show that
$\mu(\Omega_u)\rightarrow 0$ as $\|u\|\rightarrow 0$ in order to
prove Lemma \ref{l6}.

Let $\epsilon>0$. We choose $\Omega_\epsilon\subset\Omega$ a compact
subset, such that $\mu(\Omega\setminus{\Omega_\epsilon})<\epsilon$.
We denote by $\Omega_{u,\epsilon}:=\Omega_u\cap\Omega_\epsilon$. Then
it is clear that
$$C[\mu(\Omega)]^{1-p^+/q}\|u\|^{p^+}\geq\left(\int_\Omega|u|^{q}
\;dx\right)^{p^+/q}\geq\left(\int_{\Omega_{u,\epsilon}}|u|^{q}\;dx
\right)^{p^+/q}\geq\delta^{p^+}[\mu(\Omega_{u,\epsilon})]^{p^+/q}.$$
The above inequality implies that $\mu(\Omega_{u,\epsilon})
\rightarrow 0$ as $\|u\|\rightarrow 0$.

Since $\Omega_u\subset\Omega_{u,\epsilon}\cup(\Omega\setminus
\Omega_\epsilon)$ we have
$$\mu(\Omega_u)\leq\mu(\Omega_{u,\epsilon})+\epsilon$$
and $\epsilon>0$ is arbitrary. We find that $\mu(\Omega_u)\rightarrow
0$ as $\|u\|\rightarrow 0$.
This concludes the proof of Lemma \ref{l6}. \qed  

\begin{lemma}\label{l7}
The functional $J$ is coercive.
\end{lemma}
\proof
For each $u\in E$ with $\|u\|>1$ by condition (A5), relation \eq{L4} 
and inequality \eq{Hol} we have
\begin{eqnarray*}
J(u)&\geq&\int_\Omega\frac{1}{p(x)}|\nabla u|^{p(x)}\;dx-\lambda
\int_{[u>u_1]}G(x,u)\;dx-\lambda\int_{[u<u_1]}G(x,u)\;dx\\
&\geq&\frac{1}{p^+}\|u\|^{p^-}-\frac{\lambda}{\gamma}\int_{[u>u_1]}
u_1^\gamma\;dx+\frac{\lambda}{\beta}\int_{[u>u_1]}u_1^\beta\;dx-
\frac{\lambda}{\gamma}\int_{[u<u_1]}u_+^\gamma\;dx+\frac{\lambda}
{\beta}\int_{[u<u_1]}u_+^\beta\;dx\\
&\geq&\frac{1}{p^+}\|u\|^{p^-}-\frac{\lambda}{\gamma}\int_{\Omega}
u_1^\gamma\;dx-\frac{\lambda}{\gamma}\int_{\Omega}u_+^\gamma\;dx\\
&\geq&\frac{1}{p^+}\|u\|^{p^-}-\frac{\lambda}{\gamma}[\mu(\Omega)]^
{1-\gamma/p^-}C_1\|u\|^\gamma-C_2\\
&\geq&\frac{1}{p^+}\|u\|^{p^-}-C_23\|u\|^\gamma-C_2,
\end{eqnarray*}
where $C_1$, $C_2$ and $C_3$ are positive constants. Since 
$\gamma<p^-$ the above inequality implies that $J(u)\rightarrow
\infty$ as $\|u\|\rightarrow\infty$, that is, $J$ is coercive.
The proof of Lemma \ref{l7} is complete.  \qed

The following result yields a sufficient condition which ensures that a weakly convergent sequence
in $E$ converges strongly, too.

\begin{lemma}\label{l8}
Assume that the sequence $(u_n)$ converges weakly to $u$ in $E$ and
$$\limsup\limits_{n\rightarrow\infty}\int_\Omega a(x,\nabla u_n)
\cdot(\nabla u_n-\nabla u)\;dx\leq 0.$$
Then $(u_n)$ converges strongly to $u$ in $E$.
\end{lemma}

\proof
Using relation \eq{ecc4} we have that there exists a positive
constant $c_5$ such that
$$A(x,\xi)\leq c_5(|\xi|+|\xi|^{p(x)}),\;\;\;\forall x\in
\overline\Omega,\;\xi\in\RR^N.$$
The above inequality implies
\begin{equation}\label{ecc5}
A(x,\nabla u_n)\leq c_5(|\nabla u_n|+|\nabla u_n|^{p(x)}),\;\;\;
\forall x\in\overline\Omega,\;n.
\end{equation}
The fact that $u_n$ converges weakly to $u$ in $E$ implies that 
there exists $R>0$ such that $\|u_n\|\leq R$ for all $n$. By 
relation \eq{ecc5}, inequalities \eq{Hol}, \eq{L4} and \eq{L5}   
we deduce that $\{\int_\Omega A(x,\nabla u_n)\;dx\}$ is 
bounded. Then, up to to a subsequence, we deduce that
$\int_\Omega A(x,\nabla u_n)\;dx\rightarrow c$. By Lemma \ref{l3} 
we obtain
$$\int_\Omega A(x,\nabla u)\;dx
\leq\liminf\limits_{n\rightarrow\infty}\int_\Omega 
A(x,\nabla u_n)\;dx=c.$$
On the other hand, since $\Lambda$ is convex, we have
$$\int_\Omega A(x,\nabla u)\;dx\geq\int_\Omega A(x,\nabla u_n)
\;dx+\int_\Omega a(x,\nabla u_n)\cdot(\nabla u-\nabla u_n)
\;dx.$$
Next, by the hypothesis  $\limsup\limits_{n\rightarrow\infty}
\int_\Omega a(x,\nabla u_n)\cdot(\nabla u_n-\nabla u)\;dx\leq 0$, 
we conclude that $\int_\Omega A(x,\nabla u)\;dx=c$.

Taking into account that $(u_n+u)/2$ converges weakly to $u$ 
in $E$ and using Lemma \ref{l3} we have
\begin{equation}\label{ecc6}
c=\int_\Omega A(x,\nabla u)\;dx
\leq\liminf\limits_{n\rightarrow\infty}
\int_\Omega A\left(x,\nabla \frac{u_n+u}{2}\right)\;dx.
\end{equation}
We assume by contradiction that $u_n$ does not converge to $u$
in $E$. Then by \eq{L6} it follows that there exist $\epsilon>0$ and 
a subsequence $(u_{n_m})$ of $(u_n)$ such that 
\begin{equation}\label{ecc7}
\int_\Omega|\nabla(u_{n_m}-u)|^{p(x)}\;dx\geq\epsilon,\;\;\;
\forall m.
\end{equation}
By condition (A4) we have
\begin{equation}\label{ecc8}
\frac{1}{2}A(x,\nabla u)+\frac{1}{2}A(x,\nabla u_{n_m})-A\left(x,\nabla
\frac{u+u_{n_m}}{2}\right)\geq k|\nabla(u_{n_m}-u)|^{p(x)}.
\end{equation} 
Relations \eq{ecc7} and \eq{ecc8} yield
$$\frac{1}{2}\int_\Omega A(x,\nabla u)\;dx+\frac{1}{2}\int_\Omega 
A(x,\nabla u_{n_m})\;dx-\int_\Omega A\left(x,\nabla
\frac{u+u_{n_m}}{2}\right)\;\geq k\int_\Omega|\nabla(u_{n_m}-u)|^{p(x)}
\;dx\geq k\epsilon.$$
Letting $m\rightarrow\infty$ in the above inequality we obtain
$$c-k\epsilon\geq\limsup\limits_{m\rightarrow\infty}
\int_\Omega A\left(x,\nabla\frac{u+u_{n_m}}{2}\right)\;dx$$
and that is a contradiction with \eq{ecc6}. It follows that $u_n$
converges strongly to $u$ in $E$ and Lemma \ref{l8} is 
proved.  \qed

\medskip
{\sc Proof of Theorem \ref{t1} completed.}
Using Lemma \ref{l6} and the Mountain Pass Theorem (see \cite{AR}
with the variant given by Theorem 1.15 in \cite{W}) we deduce that 
there exists a sequence $(u_n)\subset E$ such that
\begin{equation}\label{MPT}
J(u_n)\rightarrow c>0\;\;\;{\rm and}\;\;\;J^{'}(u_n)\rightarrow 0
\end{equation}
where
$$c=\inf\limits_{\gamma\in\Gamma}\max\limits_{t\in[0,1]} 
J(\gamma(t))$$
and
$$\Gamma=\{\gamma\in C([0,1],E);\;\gamma(0)=0,\;\gamma(1)=u_1\}.$$
By relation \eq{MPT} and Lemma \ref{l7} we obtain that $(u_n)$ is 
bounded and thus passing eventually to a subsequence, still denoted 
by $(u_n)$, we may assume that there exists $u_2\in E$ such that 
$u_n$ converges weakly to $u_2$. Since $E$ is compactly embedded 
in $L^i(\Omega)$ for any $i\in[1,p^-]$, it follows that $u_n$ 
converges strongly to $u_2$ in $L^i(\Omega)$ for all $i\in[1,p^-]$. 
Hence
$$\langle\Lambda^{'}(u_n)-\Lambda^{'}(u_2),u_n-u_2\rangle=
\langle J^{'}(u_n)-J^{'}(u_2),u_n-u_2\rangle+\lambda\int_\Omega
[g(x,u_n)-g(x,u_2)](u_n-u_2)\;dx=o(1),$$
as $n\ri\infty$.
By Lemma \ref{l8}  we deduce that $u_n$ converges strongly
to $u_2$ in $E$ and using relation \eq{MPT} we find
$$J(u_2)=c>0\;\;\;{\rm and}\;\;\;J^{'}(u_2)=0.$$ 
Therefore, $J(u_2)=c>0$ and $J^{'}(u_2)=0$. By Lemma \ref{l5} we
deduce that $0\leq u_2\leq u_1$ in $\Omega$. Therefore
$$g(x,u_2)=u_2^{\gamma-1}-u_2^{\beta-1}\;\;\;{\rm and}\;\;\;
G(x,u_2)=\frac{1}{\gamma}u_2^\gamma-\frac{1}{\beta}u_2^\beta$$
and thus 
$$J(u_2)=I(u_2)\;\;\;{\rm and}\;\;\;J^{'}(u_2)=I^{'}(u_2).$$ 
We conclude that $u_2$ is a critical point of $I$ and thus a 
solution of \eq{2}. Furthermore, $I(u_2)=c>0$ and $I(u_2)>0>I(u_1)$. 
Thus $u_2$ is not trivial and $u_2\neq u_1$.
The proof of Theorem \ref{t1} is now complete.
\qed

\end{document}